\documentclass[12pt]{article}
\usepackage{amsmath,amssymb}
\def \1{{\bf 1}}
\def \Adm{{\rm Adm}}
\def \al{\alpha}
\def \bs{\backslash}
\def \C{{\mathbb C}}
\def \CC{{\cal C}}
\def \CE{{\cal E}}
\def \CF{{\cal F}}

\def \CO{{\cal O}}
\def \CP{{\cal P}}
\def \CV{{\cal V}}
\def \df{\ \begin{array}{c} _{\rm def}\\ ^{\displaystyle =}\end{array}\ }

\def \ga{\gamma}
\def \Ga{\Gamma}
\def \Hom{{\rm Hom}}
\def \Id{{\rm Id}}
\def \la{\lambda}
\def \l{\ell}
\def \N{{\mathbb N}}
\def \ph{\varphi}
\def \prf{{\bf Proof:\ }}
\def \qed{\ifmmode\eqno $\square$\else\noproof\vskip 12pt plus 3pt minus 9pt \fi}
 \def\noproof{{\unskip\nobreak\hfill\penalty50\hskip2em\hbox{}%
     \nobreak\hfill $\square$\parfillskip=0pt%
     \finalhyphendemerits=0\par}}
\def \ra{\rightarrow}

\def \tr{{\rm tr\,}}
\def \vol{{\rm vol}}
\def \={\ =\ }
\def \Z{{\mathbb Z}}

\newtheorem{theorem}{Theorem}[section]

\newtheorem{lemma}[theorem]{Lemma}
\newtheorem{corollary}[theorem]{Corollary}
\newtheorem{proposition}[theorem]{Proposition}

\newcommand{\norm}[1]{\left|\hspace{-1pt}\left| #1\right|\hspace{-1pt}\right|}

\begin{document}

\pagestyle{myheadings} \markright{A DISCRETE LEFSCHETZ FORMULA}

\title{A discrete Lefschetz formula}
\author{Anton Deitmar}
\date{}
\maketitle

{\bf Abstract:}
A Lefschetz formula is given that relates loops in a regular finite graph to traces of
a certain representation.
As an application the poles of the Ihara/Bass zeta function are expressed as dimensions of
global section spaces of locally constant sheaves.

\tableofcontents

\newpage
\begin{center}
{\bf Introduction}
\end{center}

In the nineteensixties Yasutaka Ihara \cite{Ihara-1, Ihara-2, Ihara-3, Ihara-4} defined a
Selberg-type zeta function for rank one $p$-adic groups.
This zeta function can also be considered as a zeta function attached to a quotient of the
Bruhat-Tits building which is a special graph. 
The theory was later extended by Ki-Ichiro Hashimoto
\cite{Hashimoto-1, Hashimoto-2, Hashimoto-3} to cover arbitrary bipartite graphs and by
Hyman Bass \cite{Bass} to arbitrary graphs. For further results see
\cite{StarkTerras-1,StarkTerras-2}.

For a finite graph $Y$ this zeta function is defined as
$$
Z_Y(u)\df\prod_p \left(1-u^{l(p)}\right)^{-1},
$$
where the product runs over all prime loops in $Y$ (see section \ref{tree}). Bass showed
that
$Z_Y(u)$ extends to a rational function without zeros, so the inverse, $Z_Y(u)^{-1}$, is a
polynomial.

In this paper we use Harmonic Analysis of the automorphism group $G$ of the universal
covering tree $X$ of $Y$ to describe the zeta function.
We prove a Lefschetz formula (Theorem \ref{5.1}) that relates lengths of loops to the
trace of co-invariants of the induced representation. This Lefschetz formula gives rise to
two descriptions of the poles of the zeta function. Firstly, the poles are given by
co-invariants of the regular $G$-representation given by the fundamental group, and
secondly they can be described in terms of a family of locally constant sheaves $\CF_\la$,
$\la\in\C$ arising from homogeneous line bundles on the boundary of the tree. The
relation is that the pole-order of $Z_Y(u)$ at the point $u=\la$ is given by
$$
{\rm ord}_{u=\la}Z_Y(u)\= \dim H^0(\CF_\la).
$$
This assertion may be viewed as a discrete version of the Patterson conjecture
\cite{BunkeOlbrich}.

The method of proof is an extension of the Harmonic Analysis of $p$-adic groups to the
automorphism group $G$.
We define the analogues of parabolic groups, split tori, and the like. It turns out that a
parabolic $P$ satisfies a Langlands decomposition $P=MAN$, but the group $M$ does not
centralize the ``split torus'' $A$. Therefore we are forced, when considering ``Jacquet
modules'', to take co-invariants not with respect to the ``unipotent radical'' $N$ but with
respect to the group $S=MN$.

There are two appendices in which generalizations of the theory of Ihara/Bass zeta
functions are suggested.
The first extends the well known generalization to twists by locally constant sheaves. It
is shown that certain constructible sheaves may serve the same purpose.
In the second the theory is extended to infinite graphs of the type that occurs as
quotients of Bruhat-Tits trees of arithmetic groups in positive characteristic
\cite{Serre}.

In writing this paper up I did not strive for the most general version of each assertion,
simply because it is much easier to read this way. Advanced readers may convince
themselves that the theory extends to multigraphs, ie, where multiple edges and
non-regular edges are allowed and that the zeta function may also be twisted by hermitian
locally constant sheaves of finite rank, ie, by finite dimensional unitary representations
of the fundamental group.
It is, however, unclear to me how far the regularity condition can be relaxed, ie, whether
the results of the paper (or modified versions thereof) hold for general finite graphs.

\section{The zeta function}
In this
paper, a \emph{graph} $X$ consists of a set $VX$, called the set of \emph{vertices},
together with a set $EX$ of subsets of $VX$ of order two, called \emph{edges}. We say that
the edge $e=\{ x,y\}$ \emph{connects} the vertices $x$ and $y$ and we also write
$e=\overline{x,y}$. We then say that $x$ and $y$ are \emph{adjacent}. The graph is called
a \emph{finite graph} if $VX$ is a finite set. For a vertex $x$ the number
$val(x)\in\{0,1,2,\dots,\infty\}$ of edges ending at $x$ is called the \emph{valency} of
$x$. The graph is called \emph{univalent} or \emph{regular} if the valency is finite and
the same for all vertices.

Sometimes we will also consider a graph as a one dimensional CW-complex, so all edges are identified with the unit interval. It does then make sense to speak of points of the graph and in particular of midpoints of edges. Further we have a path-length metric $d(.,.)$ on $X$ giving each edge the length one.

A (combinatorial) \emph{path} $p$ is a sequence $(v_0,\dots,v_n),\ n>0$ of vertices such
that
$v_j$ and
$v_{j+1}$ are adjacent for $j=0,\dots, n-1$. The number $n=l(p)$ is called the
\emph{length} of the path. Likewise, a path can also be described by the oriented edges
it passes through. The path is
\emph{closed} if
$v_0=v_n$. The path has
\emph{no backtracking} if $v_{j-1}\ne v_{j+1}$ for every $j=1,\dots,n-1$. 
 A closed path is \emph{tail-less} if $v_1\ne v_{n-1}$. A tail-less closed path with no
backtracking is called \emph{reduced}. Two reduced paths $(v_0,\dots v_n)$ and
$(w_0,\dots,w_k)$ are said to be \emph{shift-equivalent} if $k=n$ and there is an
$s=0,1,2,\dots$ such that $v_j=w_{j+k}$ where $j$ runs modulo $n$. 
 A shift equivalence class of reduced paths is called a \emph{loop}. 
 If $c$ is a loop, then $c^2=cc$ also is a loop and so is any power $c^k$,
$k\in\N$. A loop $c$ is called \emph{prime} if it is not a power $d^k$ of a
shorter loop $d$. 
 Every loop $c$ is a power $c=c_0^{m(c)}$ of a unique prime loop $c_0$. 
 The number $m(c)\in\N$ is called the \emph{multiplicity} of $c$.
 The length $l(c_0)$ is called the \emph{geometric length} of the loop $c$.
Then obviously, $l(c)=m(c)l(c_0)$. Finally, a graph is \emph{connected} if any two
vertices can be connected by a path.
 
 Let $Y$ be a finite connected graph. Its \emph{zeta function} is defined by
 $$
 Z_Y(u)\df \prod_p \left(1-u^{l(p)}\right)^{-1},
 $$
 where the product runs over all prime loops in $Y$. It has been shown by H.
Bass in \cite{Bass} that $Z_Y(u)$ always is a rational function, indeed, its
inverse $Z_Y(u)^{-1}$ is a polynomial. 
 In this paper we want to study $Z_Y(u)$ of a regular finite graph via the universal covering $X$ of $Y$.

\section{The tree}\label{tree}
Let $X$ be a tree, ie, $X$ is a countable, connected graph that admits
no loops.  We assume that $X$ is regular. Let $q+1$ be its valency where
$q\in\N$.  Let $G$ be the group of automorphisms of $X$.  
The group $G$ carries a totally disconnected topology. A basis of neighbourhoods of the unit is given by the compact open subgroups
$$
K_E \df \{ g\in G: ge=e\ \forall e\in E\}
$$
where $E$ is an arbitrary nonempty finite set of vertices. Every compact subgroup of $G$ fixes a point $x$ which is either a vertex or the midpoint of an edge. So every maximal compact subgroup is a stabilizer $K_x$ of such a point $x$. 
Another way to describe the topology of $G$ is to say that a sequence $(g_n)$ in $G$ converges to $g\in G$ iff for every vertex $x$ there is a natural number $N(x)$ such that for every $n\ge N(x)$ we have $g_n(x)=g(x)$.
 
A \emph{ray} in $X$ is a one sided infinite path without backtracking, so it is a sequence of vertices $(v_0,v_1,\dots)$ such that $v_j$ and $v_{j+1}$ are adjacent and $v_j\ne v_{j+2}$ for $j\ge 0$. A \emph{line} is a two-sided infinite path $(\dots,v_{-1},v_0,v_1,\dots )$ without backtracking.
 
Two rays $(v_0,v_1,\dots)$ and $(w_0,w_1,\dots)$ are called \emph{parallel} if there are $s\in\Z$ and $N\in\N$ such that $v_j=w_{j+s}$ for every $j\ge N$. Parallelity is an equivalence relation. The \emph{boundary} $\partial X$ of $X$ is by definition the set of all parallelity classes of rays in $X$.
A boundary point $c\in\partial X$ also is called a \emph{cusp} of $X$.

The group $G$ acts on the boundary $\partial X$. A \emph{parabolic subgroup} $P$ of $G$ is the stabilizer of a cusp $b\in\partial X_{geo}$. We fix such a pair $(P,b)$. Every vertex $x$ can be joined by a ray to $b$. In other words, given $x$, there exists a ray $(x=v_0,v_1,\dots)$ which lies in the parallelity class $b$. 
This ray is unique. We write it as $(v_0(x),v_1(x),\dots)$. 
Two vertices $x,y$ are said to lie in the same \emph{horosphere} wrt $b$ if there is $N\in\N$ such that $v_N(x)=v_N(y)$. Then automatically, $v_{N+j}(x)=v_{N+j}(y)$ for $j\ge 0$. 
Loosely speaking a horosphere is the set of all vertices of the same distance to $b$. If $p\in P$ preserves a horosphere, it will preserve every horosphere. 
Let $S\subset P$ be the subgroup of those elements that preserve horospheres.

Let  $l=(\dots, v_{-1},v_0,v_1,\dots)$ be a line prolonging the ray $(v_0,v_1,\dots)$.
We keep this line fixed and denote by $\infty\in \partial X$ the cusp given by the ray
$(v_0,v_1,v_2,\dots)$ and by $0$ the cusp corresponding to $(v_0,v_{-1},v_{-2},\dots)$. Let
$L\subset P$ be the subgroup that preserves the line $l$. It is called a \emph{Levi
component} of
$P$. Let $\bar P$ be the stabilizer of the cusp $0$. Then $L=P\cap \bar P$.

Let $M\df L\cap S$. Then $M$ is the pointwise stabilizer of the line $l$. The group $M$ is
normal in
$L$ and
$L/M$ is infinite cyclic. Therefore the exact sequence
$$
1\ \ra\ M\ \ra\ L\ \ra\ C\ \ra\ 1
$$
splits. Here $C$ denote an infinite cyclic group. Fix a splitting $s\colon C\ra L$ and let
$A$ denote the image of $s$. Then $L=AM=MA$. Note that $M$ is compact. 

Fix a finite ring $R$ with $q$ elements. At each vertex $v$ there are $q$ edges which do not point towards $\infty$. 
We label these edges with the elements of $R$ in a way that
\begin{itemize}
\item 
for $v\in l$ the edge pointing towards the cusp $0$ is labelled by $0\in
R$.
\item
The labelling is invariant under the action of $A$.
\end{itemize}

For each labelled edge $e$ let $lab(e)\in R$ denote the label.
For $b\in\partial X\setminus \{\infty\}$ there is a unique line $(\dots,w_{-1},w_0,w_1\dots)$ joining $b$ to $\infty$.
Define $\psi(b)$ in the ring of formal Laurent series $R((t))$ by
$$
\psi(b)\df \sum_{j=-\infty}^\infty lab(\overline{w_j,w_{j-1}})\, t^j.
$$
The map $\psi$ is a bijection from $\partial X\setminus \{\infty\}$ to the ring $R((t))$.

Note that the action of $S$ on $\partial X\setminus \{\infty\}$ defines an injection $S\hookrightarrow {\rm Per}(\partial X\setminus \{\infty\})$. Let $N$ denote the subgroup of $S$ consisting of all $n\in S$ for which there is $n_R\in R((t))$ such that
$$
n(b)\= \psi^{-1}(\psi(b)+n_R).
$$
In other words, $N$ consists of those elements of $S$ which act by addition of an element of $R((t))$ on $\partial X$.
The group $N$ thus is isomorphic to the additive group of $R((t))$ and acts transitively on $\partial X\setminus \{\infty\}$. 
It follows that $S=MN$ and $P=MAN$.
The group $N$ is a closed subgroup of $G$. Indeed, it carries the $t$-adic topology of
$R((t))$.

Let $\bar P$ denote the stabilizer of the cusp $0$.
Define $\bar N$ in an analogous fashion to $N$ to get $\bar P=MA\bar N$.

There is an element $w\in G$ with $w^2=1$ and $wv_j=v_{-j}$ for every $j$. Further $w$ is
supposed to be chosen so that
$$
waw\= waw^{-1}\= a^{-1}
$$
for every $a\in A$. Then $w$ is called \emph{the non-trivial Weyl element} and $W=\{
1,w\}$ is the \emph{Weyl group} to $A$.

\begin{lemma}
(Bruhat decomposition)\\
The group $G$ can be decomposed into disjoint sets,
$$
G=P\cup PwP.
$$
Further, $PwP=NwP=PwN$.
\end{lemma}

\prf
The group $N$ acts transitively on the cusps different from $\infty$. This implies that
for every $g\in G\setminus P$ there is $n\in N$ with $ng\infty=0$.
This implies $wng\in P$ and so $g\in NwP$.
\qed

\section{Trace formula and duality}
Let $G$ be a totally disconnected locally compact group, so $G$ has a basis of unit neighbourhoods consisting of compact open subgroups.
By a \emph{representation} of $G$ we mean a group homomorphism from $\pi\colon G\ra {\rm GL}(V_\pi)$ for some complex vector space $V_\pi$. By abuse of notation we will write $\pi$ for the space $V_\pi$ as well. 
For every compact open subgroup $K$ of $G$ let $\pi^K$ denote the space of $K$-fixed vectors.
We say that the representation $\pi$ is \emph{smooth} if every $v\in\pi$ is
fixed by some compact open subgroup $K$. 

Let $\pi$ be a smooth representation. For a compact subgroup $K$ of $G$ and
$v\in \pi$ let
$$
\CP_K(v)\df\frac 1{\vol(K)}\int_K\pi(k)v\, dk,
$$
the projection to the $K$-invariants.
In other words, $\CP_K$ is a linear map on $\pi$ given as follows.
For given $v\in\pi$ there exists a subgroup $K'$ of finite index in $K$ such
that $v$ is invariant under $K'$. Then
$$
\CP_K(v)\=\frac 1{|K/K'|}\sum_{x\in K/K'}\pi(x)v.
$$
It is easy to see that $\pi$ is the direct sum of the space $\pi^K$ and the
kernel $\pi[K]$ of the linear map $\CP_K$.

A smooth representation $\pi$ is
\emph{admissible} if for every compact open subgroup $K$ the space $\pi^K$ is
finite dimensional.

Suppose that $G$ admits  a uniform lattice $\Ga$, ie, a discrete subgroup such that $\Ga\bs G$ is compact. 
Then $G$ is unimodular (Ex 14.2 in \cite{Dieu}). 
The space of $C^\infty(\Ga\bs G)$ locally constant functions on $\Ga\bs G$
carries an admissible representation $R$ given by
$$
R(y)\ph(x)\df \ph(xy).
$$
 
 A function $f$ on $G$ is called a \emph{Hecke function} if $f$ is integrable
and there is a compact open subgroup $K$ of $G$ such that $f$ factors over
$K\bs G/K$. For any Hecke function $f$ the integral
$$
R(f)\df\int_G f(x)\, R(x)\, dx
$$
defines a  linear operator on $C^\infty(\Ga\bs G)$. For $g\in G$ and $f$ a
function on $G$ define the \emph{orbital integral}
$$
\CO_g(f)\df \int_{G_g\bs G}f(x^{-1}gx)\, dx
$$
whenever the integral exists. Here $G_g$ is the centralizer of $g$ in $G$. We have to fix a Haar measure on $G_g$ here.

\begin{proposition}
(Trace formula)\\
Let $f$ be a Hecke function. Then
$$
 \tr R(f)\= \sum_{[\ga]} \vol(\Ga_\ga\bs
G_\ga)\, \CO_\ga(f),
$$
where the sum on the right hand side runs over the set of all conjugacy classes in the
group $\Ga$ and $\Ga_\ga$ denotes the centralizer of $\ga$ in
$\Ga$.
\end{proposition}

The proof is the same as that of the corresponding assertion for $p$-adic
linear groups \cite{padgeom}.
\qed

The trace formula can also be used to show the Ihara/Bass-identity (see \cite{Terras,
Venkov-Nikitin}).

We now turn to duality. Let $\Adm(G)$ denote the category of all admissible smooth
representations of $G$.
We introduce the setup of topological (continuous) representations.
Let $\CC(G)$ denote the category of continuous representations of $G$ on locally convex, complete, Hausdorff topological vector spaces.
The morphisms in $\CC(G)$ are continuous linear $G$-maps.

To any $V\in \CC(G)$ we can form $V^\infty$, the space of smooth vectors, which is by
definition the space of all vectors in $V$ which have an open stabilizer. By the
continuity of the representation it follows that $V^\infty$ is dense in $V$. 
We call the representation $V$ \emph{admissible} if $V^\infty$ is admissible.
Let $\Adm_{top}(G)$ denote the category of all admissible topological representations.
The
 functor:
\begin{eqnarray*}
F : \Adm_{top}(G) &\ra & \Adm(G)\\
W&\mapsto & W^\infty 
\end{eqnarray*}
is easily seen to be exact.

Let $V$ be in $\Adm(G)$ then any $W\in\Adm_{top}(G)$ with $FW=V$ will be called a
\emph {completion} of $V$.

To $V\in \Adm(G)$ let $V^*=Hom_\C(V,\C)$ be the dual module and let $\tilde{V}$ be the
space of smooth vectors in $V^*$. Then $\tilde V$ is admissible again and the
natural map from
$V$ to
$\tilde{\tilde{V}}$ is an isomorphism.

On the space $C(G)$ of continuous maps from $G$ to the complex numbers we have two actions of $G$, the left and the right action given by
$$
L(g)\ph(x):=\ph(g^{-1}x)\ \ \  R(g)\ph(x) := \ph(xg),
$$
where $g,x\in G$ and $\ph\in C(G)$. 
For $V\in\Adm(G)$ let
$$
V^{-\infty} := Hom_G(\tilde{V} ,C(G)),
$$
where we take $G$-homomorphisms with respect to the right action, so $V^{-\infty}$ is the
space of all linear maps $f$ from $\tilde{V}$ to $C(G)$ such that
$f(g.v^*)=R(g)f(v^*)$. Then $V^{-\infty}$ becomes a $G$-module via the left translation:
for $\alpha\in V^{-\infty}$ we define $g.\alpha (v^*) := L(g)\alpha(v^*)$. We call
$V^{-\infty}$ the \emph{maximal completion} of $V$. The next lemma and the next
Proposition will justify this terminology. First observe that the topology of locally
uniform convergence may be installed on $V^{-\infty}$ to make it an element of $\CC(G)$.
By $V\mapsto V^{-\infty}$ we then get a functor $R:\Adm(G)\ra\Adm_{top}(G)$.

\begin{lemma}
We have $(V^{-\infty})^\infty \cong V$.
\end{lemma}

\prf
There is a natural map $\Phi : V\ra V^{-\infty}$ given by $\Phi(v)(v^*)(g) := v^*(g^{-1}.v)$, for $v^*\in \tilde{V}$ and $g\in G$.
This map is clearly injective.
To check surjectivity let $\alpha\in(V^{-\infty})^\infty$, then the map
$\tilde{V}\ra \C$, $v^*\mapsto \alpha(v^*)(1)$ lies in $\tilde{\tilde{V}}$.
Since the natural map $V\ra \tilde{\tilde{V}}$ is an isomorphism, there is a $v\in V$ such that $\alpha(v^*)(1)=v^*(v)$ for any $v^*\in \tilde{V}$, hence $\alpha(v^*)(g)= \alpha(g.v^*)(1)=g.v^*(v)=v^*(g^{-1}.v)=\Phi(v)(v^*)(g)$, hence $\alpha =\Phi(v)$.
\qed

It follows that the functor $R$ maps $\Adm(G)$ to $\Adm_{top}(G)$ and that $FR=Id$.

\begin{proposition}
The functor $R: \Adm(G)\ra \Adm_{top}(G)$ is right adjoint to $F: W\mapsto W^\infty$.
So for $W\in \Adm_{top}(G)$ and $V\in \Adm(G)$ we have a functorial isomorphism:
$$
Hom_{G}(FW,V)\cong Hom_{\CC(G)}(W,RV).
$$
\end{proposition}

\prf
We have a natural map
\begin{eqnarray*}
Hom_{\CC(G)}(W,V^{-\infty}) &\ra& Hom_{G}(W^\infty ,V)\\
\alpha &\mapsto& \alpha |_{W^\infty}.
\end{eqnarray*}
Since $W^\infty$ is dense in $W$ this map is injective.
For surjectivity we first construct a map $\psi : W\ra (W^\infty)^{-\infty}$.
For this let $W'$ be the topological dual of $W$ and for $w'\in W'$ and $w\in W$ let
$$
\psi_{w',w}(x) := w'(x^{-1}.w),\ \ \ x\in G.
$$
The map $w\mapsto \psi_{.,w}$ gives an injection
$$
W\hookrightarrow Hom_G(W',C(G)).
$$
Let $\ph\in \tilde{FW}$, then $\ph$ factors over $(FW)^K$ for some compact open subgroup $K\subset G$.
But since $(FW)^K =W^K$ it follows that $\ph$ extends to $W$.
The admissibility implies that $\ph$ is continuous there.
So we get $\tilde{(FW)} \hookrightarrow W'$.
The restriction then gives
$$
Hom_G(W',C(G))\ra Hom_G(\tilde{(FW)},C(G)).
$$
From this we get an injection
$$
\psi : W\hookrightarrow (FW)^{-\infty},
$$
which is continuous.
Let $\zeta : W^\infty \ra V$ be a morphism in $\Adm(G)$.
We get
$$
\alpha(\zeta) : W\hookrightarrow (FW)^{-\infty} \begin{array}{c} \zeta^{-\infty}\\ \longrightarrow\\ {}\end{array} V^{-\infty},
$$
with $\alpha(\zeta)|_{W^\infty}=\zeta$, i.e. the desired surjectivity.
\qed

Let now $\Ga$ denote a cocompact torsion-free discrete subgroup of $G$ and let $C(\Ga \bs
G)$ denote the space of all continuous functions on $\Ga\bs G$.

The next theorem is a version of the classical duality theorem
\cite{Gelfand} adapted to our context.

\begin{theorem}\label{higher_duality}
(Duality Theorem) 
Let $\Ga$ be a cocompact torsion-free discrete subgroup of $G$, then for any $V\in
\Adm(G)$:
$$
H^0(\Ga ,V^{-\infty}) 
	\ \cong\ \Hom_{G}(C^\infty(\Ga\bs G) ,V).
$$
\end{theorem}

There is a version for higher cohomologies in \cite{padgeom} but this won't be needed here.

\prf
Since $V^{-\infty}=\Hom_G(\tilde V,C(G))$ it follows that
\begin{eqnarray*}
H^0(\Ga,V^{-\infty}) &=& \Hom_G(\tilde V,C(\Ga\bs G))\\
&=& \Hom_G(\tilde V,C^\infty(\Ga\bs G)),
\end{eqnarray*}
since any $G$-homomorphism from the smooth $\tilde V$ to $C(\Ga\bs G)$ has image in
$C^\infty(\Ga\bs G)$.
Dualizing gives
$$
\dim H^0(\Ga,V^{-\infty})\=\dim \Hom_G(C^\infty(\Ga\bs G),V).
$$
\qed

\section{A character identity}
The proofs in this section are inspired by \cite{Casselman-1, casselman-2}.

Fix a line $(\dots,v_{-1},v_{0},v_{1},\dots)$. Let $\infty$ denote the cusp given by the ray $(v_{0},v_1,\dots)$ and let $0$ denote the cusp given by the ray $(v_0,v_{-1},v_{-2},\dots)$. 
Let $P,L,M,S,A$ be the groups defined above.

For $j\ge 0$ let $K_j$ be the group of all $g\in G$ with $gv=v$ for every vertex $v$ that satisfies $d(v_0,v)\le j$.
Then $K_0$ is a maximal compact subgroup and the $K_j$ form a basis of neighbourhoods of the unit in $G$.
Let $M_j=M\cap K_j$, $N_j=N\cap K_j$ and $\bar N_j=\bar N\cap K_j$.

For any set $X$ let $\1_X$ denote its characteristic function.
For $f,g\in L^1(G)$ denote by $f*g$ their convolution product, so
$$
f*g(x)\= \int_G f(y)g(y^{-1}x)\, dy.
$$

\begin{lemma}\label{4.1}
Let $j\ge 1$.
\begin{enumerate}
\item 
We have the Iwahori factorization,
$K_j=\bar N_j M_j N_j= N_jM_j\bar N_j$.
\item
For $a\in A^-$,
\begin{itemize}
\item 
$aM_ja^{-1}=M_j$,
\item
$a N_j a^{-1}= N_{j+l(a)}$,
\item
$ a^{-1}\bar N_j a = \bar N_{j+l(a)}$.
\end{itemize}
\item
For $a\in A$,
$$
\vol(K_jaK_j)\= q^{l(a)}\, \vol(K_j).
$$
\item
For $a,b\in A^-$,
$$
\1_{KaK}*\1_{KbK}\= \1_{KabK}.
$$
\end{enumerate}

\end{lemma}

\prf
For (a) note that if $k\in K_j,\ j\ge 1$, then $k0\ne \infty$.
hence there are $n_j\in N_j$ and $\bar n_j\in \bar N_j$ such that
\begin{eqnarray*}
\bar n_j k n_j \infty &=& \infty,\\
\bar n_j k n_j 0 &=& 0,
\end{eqnarray*}

which implies $\bar n_j k n_j \in M_j$. Now the first part of (a) follows. The second is obtained by taking
inverses. The point (b) is immediate from the fact that $av_j=v_{j+v(a)}$ and $N_j$ being
the stabilizer of $v_j$ in $N$.

For (c) note that the map $k_1 ak_2\mapsto k_1$ induces a bijection between the sets $K_jaK_j/K_j$ and $K_j/(K_j\cap aK_j a^{-1})$.
But since $a\in A^-$, $K_j\cap aK_ja^{-1}=\bar N_j M_j aN_ja^{-1}=\bar N_j M_j N_{j+l(a)}$ and $[K_j:K_j\cap aKa^{-1}]=[N_j:N_{j+l(a)}]=q^{l(a)}$.

Finally, we show (d).
As sets one has $K_jaK_j\cdot K_jbK_j=K_j ab K_j$ because for $k_j=n_jm_j\bar n_j$ one has $ ak_jb=an_ja^{-1}\cdot ab\cdot b^{-1}m_jb\cdot b^{-1}n_jb$.
By (c) above also the measures agree.
\qed

Let $\pi$ be an admissible representation of $G$.
For any subgroup $H$ on $G$ 
let $\pi[H]$ denote the complex vector space generated by all elements of the form $v-\pi(h)v$, $v\in \pi$, $h\in H$. 
For $H=S$ let $\pi_S$ be the quotient $\pi /\pi[S]$. Then $\pi_S$ is the largest $P$-module quotient on $\pi$ on which $S$ acts trivially.

\begin{corollary}\label{4.2}
For $v\in \pi^{K_j}$ with image $u$ in $\pi_S$ the vector 
$$
\frac 1{\vol(K_jaK_j)} \pi(\1_{K_jaK_j})v
$$ 
has image $\pi_S(a)u$.
\end{corollary}

\prf 
This follows from Lemma \ref{4.1}.
\qed

\begin{proposition}
Suppose $\pi$ is finitely generated as $G$-module. Then $\pi_S$ is a finitely generated $A$-module.
\end{proposition}

\prf
Let $E$ be a finite generating set of $\pi$.
Let $K\subset G$ be a compact open subgroup such that $E\subset \pi^K$.
Let $\Lambda$ be a finite subset of $G$ such that $P\Lambda K =G$.
Since $\pi$ is the linear span of $\pi(G)E$ it follows that $\pi_S$ is as $A$-module generated by the image of $\pi(\Lambda) E$.
\qed

\begin{proposition}\label{4.4}
\begin{enumerate}
\item 
If $v\in\pi$ is fixed by $M_j\bar N_j$, then $\CP_{K_j}(v)=\CP_{N_j}(v)$.
\item 
The canonical projection from $\pi^M$ to $\pi_S$ is surjective.
\item 
The canonical projection from $\pi^{K_1}$ to $\pi_S$ is surjective. In particular, $\pi_S$ is finite dimensional.
\end{enumerate}

\end{proposition}

\prf
The first assertion is an easy consequence of the Iwahori decomposition.
Let $\pi_M$ be the largest quotient on which $M$ acts trivially. 
Since $M$ is compact, integration over $M$ induces an isomorphism $\pi_M\ra\pi^M$. Since the map $\pi_M\ra\pi_S$ is surjective, the second claim  follows.

For the third let $U_S$ be a finite dimensional subspace of $\pi_S$ and let $U$ be a finite dimensional subspace of $\pi^M$ mapping to $U_S$. 
There is $j\in\N$ such that $U\subset \pi^{M_1\bar N_j}$.
Choose $a\in A$ such that $a^{-1}\bar N_1 a\subset \bar N_j$.
Then $\pi(a) U\subset \pi^{M_1\bar N_1}$.
The decomposition $K_j=\bar N_j M_j N_j$ for $j=1$ implies that the image of $ \pi^{M_1\bar N_1}$ in $\pi_S$ equals the image of $\pi^{K_1}$.
So it follows that $\pi(a) U_S$ lies in the image of $\pi^{K_1}$.
Hence its dimension, ie, the dimension of $U_S$, is bounded by the dimension of $\pi^{K_1}$.
We conclude that $\pi_S$ is finite dimensional.
Thus one can take $U_S$ to be equal to $\pi_S$.
But then $U_S=\pi(a)U_S$ lies in the image of $\pi^{K_1}$.
\qed

\begin{proposition}\label{4.5}
Let $a_1$ be the unique generator of $A$ inside $A^-$. For every $j\ge 1$ there exists
a space $\pi_1^{K_j}\subset \pi^{K_j}$ such that
\begin{enumerate}
\item 
The projection from $\pi_1^{K_j}$ to $\pi_S$ is a linear isomorphism.
\item
For each $m\ge 0$ the space $\pi_1^{K_j}$ is stable under
$\pi(\1_{K_j a_1^m K_j})$.
\item
There exists $n\in\N$ such that $\pi(\1_{K_j a_1^nK_j})\pi^{K_j}\subset \pi_1^{K_j}$.
\end{enumerate}

\end{proposition}

\prf
We need the following lemma.

\begin{lemma}
For every compact subgroup $H$ of $G$ the space $\pi[H]$ coincides with the space of all $v\in\pi$ with $\int_H\pi(h)v\, dh =0$.
\end{lemma}

\prf
For every $v\in\pi[H]$ we obviously have $\int_H\pi(h)v\, dh=0$. The other way round assume that $v\in\pi$ satisfies $\int_H\pi(h)v\, dh=0$.
There exists a normal subgroup $H'$ of finite index in $H$ such that $v\in\pi^{H'}$.
Let $F$ denote the finite group $H/H'$. Then one has $0=\sum_{x\in F}\pi(x)v$.
Hence
$$
v\=\frac 1{|F|} \sum_{x\in F} (v-\pi(x)v).
$$
This shows that $v$ lies in $\pi[H]$.
\qed

Note that $S$ is the union of the compact open subgroups $S_j=K_j\cap S=M_jN_j$.
Thus $\pi[S]$ is the union of the sets $\pi[U]$, when $U$ runs over all compact open subgroups of $S$.
Fix $j\in\N$.
Choose a fixed compact open subgroup $U$ of $S$ such that $\pi[S]\cap\pi^{K_j}\subset \pi[U]$ and $S_j\subset U$.

\begin{lemma}\label{4.7}
If $a_1^n Ua_1^{-n} \subset S_j$ and $v\in\pi[S]\cap \pi^{K_j}$, then $\pi\left(\1_{K_j a_1^n K_j}\right) v=0$.
\end{lemma}

\prf
The vector $\pi\left(\1_{K_ka_1^nK_j}\right) v$ differs from $\CP_{K_j}(\pi(a_1^n)v)$ only by a scalar.
By Proposition \ref{4.4} the latter equals $\CP_{N_j}(\pi(a_1^n)v)$.
But
\begin{eqnarray*}
\CP_{N_j}(\pi(a_1^n)v) &=& (const) \int_{N_j}\pi(n)\pi(a_1^n)v\, dn\\
&=& (const)\, \pi(a_1^n)\int_{a_1^nN_ja_1^{-n}}\pi(n)v\, dn\= 0.
\end{eqnarray*}
\qed

Choose $n$ to be large enough that $a_1^nUa_1^{-n}\subset S_j$ and
define $\pi_1^{K_j}$ to be $\pi(K_j a_1^n K_j) V^{K_j}$.
We first show that thw projection from $\pi_1^{K_j}$ to $\pi_S$ is
surjective. For this let $u\in\pi_S$. By Proposition \ref{4.4} there is
$v\in\pi^{K_j}$ whose image in $\pi_S$ equals $\pi_S(a_1^{-n})u$.
By Corollary \ref{4.2}, $\CP_{K_j}(\pi(a_1^n)v)$ has image $u$.

For the injectivity let $v_0\in\pi^{K_j}$ and assume that $v=\pi(K_j a_1^n
K_j)v_0$ lies in $\pi[S]$. We have to show that $v=0$.

By choice of $U$, $v\in\pi[U]$. Now $v$ is also, up to a constant, equal to
$\CP_{K_j}(\pi(a_1^n)v_0)=\CP_{N_j}(\pi(a_1^n)v_0)$. Therefore
\begin{eqnarray*}
\int_U \pi(u) v\, du &=& 0\\
&=& \int_U\pi(u)\, du\, \int_{N_j}\pi(n_j)\,\pi(a_1^n) v_0\,dn_j\\
&=& \int_U \pi(u) \pi(a_1^n)v\, du\\
&=& \pi(a_1^n)\int_{a_1^{-n}Ua_1^n} \pi(u) v\, du.
\end{eqnarray*}
So $v_0\in\pi[S]$ and by Lemma \ref{4.7} it follows $v=0$.

To prove (b) note that the above is valid for every large enough $n$, so all
the spaces $\pi(K_j a_1^nK_j)\pi^{K_j}$ have the same dimension. But for
$m\ge 0$,
$$
\pi(K_j a_1^m K_j)\pi(K_ja_1^nK_j)\pi^{K_j}\=\pi(K_ja^{M+n}K_j)\pi^{K_j}.
$$
Finally, statement (c) follows from the definition.
\qed

We now state the main theorem of this section.

\begin{theorem}
For all $j,m\ge 1$,
$$
\tr\left(\frac 1{\vol(K_ja_1^mK_j)}\pi(\1_{K_j a_1^mK_j})\right) \=
\tr\pi_S(a_1^m).
$$
\end{theorem}

\prf
Since $(\1_{K_j a_1^mK_j})^n=\1_{K_ja_1^{mn}K_j}$, the Proposition \ref{4.5}
implies that 
\begin{eqnarray*}
\tr\left(\frac 1{\vol(K_ja_1^mK_j)}\pi(\1_{K_j a_1^mK_j})\right) &=& 
\frac 1{\vol(K_ja_1^mK_j)}\tr\left(\pi(\1_{K_j a_1^mK_j})\mid
\pi_1^{K_j}\right)\\
&=& \tr(a_1^m|\pi_S).
\end{eqnarray*}
\qed

An element $g$ of $G$ is called \emph{elliptic} if it fixes a point of $X$. This point can
be chosen to be a vertex or a midpoint of an edge. If $g$ fixes a midpoint, then $g^2$
fixes both vertices of that particular edge. If $g$ does not fix a point of $X$ then it is
called \emph{hyperbolic}. In that case the minimum
$$
l(g)\df \min_{x\in X} d(gx,x)
$$
is attained exactly on a line $l(g)=(\dots,v_{-1}(g),v_0(g),v_1(g),\dots)$ which is
preserved by
$g$. Since this line $l(g)$ can be mapped to our given line by some $h\in G$, it follows
that $g$ is conjugate to an element of $L=AM$.
It is not hard to see that every hyperbolic element of $L$ is conjugate to a uniquely
determined element of
$A^-$. Let
$G_{hyp}$ denote the set of hyperbolic elements of $G$.

\begin{corollary}
Let $\pi$ be an admissible representation. 
On $G_{hyp}$ define a conjugation invariant function $\theta_\pi$ by
$\theta_\pi(a)=\tr(a|\pi_S)$, $a\in A^-$.
Then for every Hecke function $f$ which is supported on $G_{hyp}$,
$$
\tr\pi(f)\=\int_Gf(x)\,\theta_\pi(x)\, dx.
$$
\end{corollary}

\prf
The claim follows from the theorem for Hecke functions of the form $f=\1_{K_j
aK_j}$, $a\in A^-$.
Every Hecke function with support in $G_{hyp}$ is a linear combination of
functions which are conjugate to functions of the latter form.
\qed

\section{The Lefschetz formula}
Recall that $G$ is the automorphism group of a $q+1$-regular tree.
Note that such a group $G$ admits a lattice $\Ga$ as above, hence $G$ is unimodular.
Let $g$ be a hyperbolic element of $G$.
Let $l(g)=l=(\dots,v_{-1},v_0,v_1,\dots)$ be the line that is preserved by $g$. After
renumbering we can achieve that $gv_j=v_{j-l(g)}$. Let $b$ be the parallelity class of the
ray
$(v_0,v_1,\dots)$. Let $P,S,M,A$ be defined as before. Then $g\in L=AM$. The centralizer
$G_g$ of $g$ preserves $l$ and therefore is a subgroup $G_g=L_g$ of $L$. Let $M_g=G_g\cap
M$, then we have an exact sequence
$$
1\ \ra\ M_g\ \ra\ L_g\ \ra\ C\ \ra\ 1
$$
for an infinite cyclic group $C$. Let $A_g\subset L_g$ be an image of a splitting $s\colon C\ra L_g$, then $L_g=A_g M_g$.
We normalize the Haar measure on the compact group $M_g$ to have volume one and on $A_g$ we choose $l(a_0)$ times the counting measure, where $a_0$ is a generator of $A_g$. Then the volume of $<g>\bs L_g$ equals $l(g)$. 

Now suppose that $\Ga$ acts freely on $X$. Then every element $\ga\in\Ga$, $\ga\ne 1$ is hyperbolic and the centralizer $\Ga_\ga$ is cyclic. Let $\ga_0$ be the unique generator of $\Ga_\ga$ with $\ga=\ga_0^{\mu}$ with $\mu>0$. Then
$$
\vol\left( \Ga_\ga\bs G_\ga\right)\= l(\ga_0).
$$
Alternatively we call an element $\tau$ of $\Ga$ a \emph{primitive} element if $\tau=\sigma^n$ for $\sigma\in\Ga$ and $n\in\N$ implies $n=1$. Then $\ga_0$ is the unique primitive element such that $\ga=\ga_0^{\mu}$ with $\mu>0$.

\begin{proposition}
Assume that $\Ga$ acts freely on $X$. Then for every Hecke function $f$,
$$
\tr R(f)\= \sum_{[\ga]} l(\ga_0)\,
\CO_\ga(f).
$$
\end{proposition}

\prf This follows from the trace formula.\qed

Fix a  cusp $b=(v_0,v_1,\dots)$ and a  line $l=(\dots,v_{-1},v_0,v_1,\dots)$
prolonging $b$.  Let $L,P,S,M,A$ be defined as before. Let $L^-$ denote the
set of all $g\in L$ such that $l(g)>0$ and $gv_j=v_{j-l(g)}$, ie, $L^-$ is
the set of elements of $L$ that ``move away'' from $b$. Let $A^-=A\cap L$.
Then $L^-=A^-M$. Let $\CE_p(\Ga)$ denote the set of all conjugacy classes
$[\ga]$ in $\Ga$ such that $\ga$ is in $G$ conjugate to an element $a_\ga
m_\ga$ of $L^-=A^- M$.

\begin{theorem}\label{5.1}(Lefschetz Formula)\\
Let $\Ga\subset G$ be a discrete cocompact subgroup which acts freely on $X$.
Then for every function $\ph$ on $A^-$ such that $\sum_{a\in A^-}\ph(a)
q^{l(a)}<\infty$,
$$
\sum_{a\in
A^-}\ph(a)\,q^{l(a)}\,\tr(a|R_S)\=\sum_{[\ga]\in\CE_b(\Ga)}l(\ga_0)\,\ph(a_\ga).
$$
\end{theorem}

The proof will occupy the rest of the section.
We start with the simple observation that the length map $l\colon G\ra \N_0$ is continuous and hence $l^{-1}(n)$ is open for every $n\in\N_0$. An element $g\in G$ is elliptic iff $l(g)=0$ and hyperbolic iff $l(g)>0$.

We next need an integration formula. Fix the Haar measure on $G$ that gives the stabilizer $K_x$ of a vertex $x$ the volume one.
On $A$ install the counting measure as Haar measure.

For $j\in\Z$ let $K(j)$ be the stabilizer of the point $v_j$ and set
$K=K(0)=K_0$. Define $N(j)=N\cap K(j)$. Then $N(j)$ is a compact open subgroup
of
$N$ that stabilizes $v_k$ for every $k\ge j$. The group $N$ is the union of
all the subgroups $N(j)$. We give $N$ the Haar measure such that $N(0)$ gets
volume one. Then $\vol(N(j))=q^j$.

\begin{lemma}
With the above normalizations we have for every $f\in L^1(G)$,
$$
\int_Gf(x)\, dx\=\int_K\int_N\int_A f(kna)\, da\,dn\,dk.
$$
If $f$ is
supported in the open set of hyperbolic elements,
$$
\int_G f(x)\, dx\= \int_{K}\int_N\int_{A^-} f(kn\, a(kn)^{-1})\, q^{l(a)}\,
da\, dn.
$$
For later use we also note that for $a\in A$, $f\in L^1(S)$,
$$
\int_S f(asa^{-1})\, ds\= q^{-v(a)}\int_S f(s)\, ds.
$$
Here $v(a)\in\Z$ is defined by $a(v_j)=v_{j+v(a)}$. Then $a\in A^-$ is
equivalent to $v(a)<0$.
\end{lemma}

\prf
On $K\bs G$ there is a unique $G$-invariant measure such that
$$
\int_Gf(x)\, dx\= \int_{K\bs G}\int_K f(kx)\, dk\, dx.
$$
The natural projection identifies $NA$ with $K\bs G$, so this measure gives a
Haar measure on the group $NA$. But $dn\, da$ also is a Haar measure, so the
uniqueness of Haar measures gives the first assertion up to a scalar.
Plugging in the test function $f=\1_K$ gives the claim.

Let $G_{hyp}$ denote the set of all hyperbolic elements
of $G$. Note that the map
\begin{eqnarray*}
G/A\times A^- &\ra& G_{hyp}\\
(x,a) &\mapsto & xax^{-1}
\end{eqnarray*}
is surjective and proper. Hence there is a nowhere vanishing function $B$ on
$G/A\times A^-$ such that
$$
\int_G f(x)\, dx\= \int_{G/A}\int_{A^-} f(xax^{-1})\, B(x,a)\, da\, dx.
$$
Since $G$ is unimodular it follows that $\int_G f(yxy^{-1})\, dx=\int_G f(x)\,
dx$ and hence $B$ does not depend on
$x$. Write $B(a)$ instead.

The projection $G\ra G/A$ identifies $KN$ with the quotient $G/A$. 
Via the first integral equality we get an identification of the measure on
$G/A$ with the measure $dk\, dn$.
Thus,
$$
\int_Gf(x)\, dx\= \int_K\int_N\int_{A^-} f(kn\, a\, (kn)^{-1})\,
B(a)\,da\,dn\,dk.
$$
We are going to compute $B(a)$ by plugging in special test functions.
First let $a_1\in A^-$ be the unique element with $l(a_1)=1$.
For $k\in \N$ write $B(k)$ instead of $B(a_1^k)$.
Fix some $a=a_1^{k}\in A^-$.
Let $f$ be the indicator function of the set $KaK\cap G_{hyp}$.
This function is $K$-central.
The set $KaK$ equals the set of all $g\in G$ that satisfy $d(gv_0,v_0)=k$.
Since there are $(q+1)q^{k-1}$ vertices in distance $k$ to $v_0$ it follows that
$$
\left| KaK/K\right|\= (q+1)q^{k-1}.
$$
An element $g$ of $KaK$ is elliptic iff $g$ fixes a point $x$ half-way between $v_0$ and $g(v_0)$. 
Since for every point $x$ in distance $k/2$ there are $q^{[k/2]}$ points in distance $k$ that lie beyond $x$ it follows that the volume of the set of elliptic elements in $KaK$ equals $q^{[k/2]}$.
Hence
$$
\int_G f(x)\, dx\= (q+1)q^{k-1}-q^{\left[\frac k2\right]}.
$$
For $a\in A$ recall that $v(a)\in\Z$ is defined by $a(v_j)=v_{j+v(a)}$. It is easy to
see that
$$
\{ [a^{-1},s] : s\in S_j\}\= S_k,
$$
where $k=\max(j,j-v(a))$.
For $a\in A^-$ and $s\in S$ we have $d(sas^{-1} v_0,v_0)=d(as^{-1}v_0,s^{-1}v_0)$ and this number equals $l(a)$ if $s\in S_0$.
If $s\in S_{j+1}\setminus S_j$ for some $j\ge 0$, then
$$
d(as^{-1}v_0,s^{-1}v_0)\= 2(j+1)+l(a).
$$
Using $\vol(S_0)=1$ and $\vol(S_{j+1}\setminus S_j)=q^j(q-1)$ we get
\begin{eqnarray*}
(q+1)q^{k-1}-q^{\left[\frac k2\right]} &=& \int_G f(x)\, dx\\
&=& B(k)+\sum_{\nu=1}^{\left[\frac{k-1}2\right]} B(k-2\nu)(q-1)q^{\nu-1}.
\end{eqnarray*}
The unique solution to this recurrence formula is $B(k)=q^k$. This implies the
second assertion of the lemma. For the third note that $\int_S f(asa^{-1})\, ds$ is a Haar
measure on $S$, therefore equals  $c(a)\int_Sf(s)\, ds$ for a constant $c(a)$. Plugging in
the test function $f=\1_K$ gives $c(a)=q^{v(a)}$.
\qed

Let $\eta\ge 0$ be compactly supported function on $N$, left invariant under
$K\cap N$ and such that $\int_N\eta(n)\, dn=1$.
Define a function $f$ on $G$ by
$$
f(kn\,a\,(kn)^{-1})\df \eta(n)\,\ph(a),\ \ \ \ a\in A^-,
$$
and $f(x)=0$ if $x$ is elliptic.
The $f$ is a Hecke function.
Let $am\in A^-M$. We have seen that $am$ is conjugate to $a$. so that
$$
\CO_{am}(f)\=\CO_a(f)\=\ph(a).
$$
If we plug this function $f$ into the trace formula we readily see that the
geometric side will give the geometric side of the Lefschetz formula.
For the spectral side we compute for an admissible representation $\pi$,
\begin{eqnarray*}
\tr\pi(f)&=&\int_Gf(x)\theta_\pi(x)\, dx\\
&=&
\int_K\int_N\int_{A^-}\eta(n)\ph(a)q^{l(a)}\tr(a|\pi_S)\,da\,dn\,dk\\
&=&\int_{A^-} \ph(a)q^{l(a)}\tr(a|\pi_S)\,da.
\end{eqnarray*}
So Theorem \ref{5.1} follows from the trace formula.
\qed

\section{Rationality of the zeta function}
Let $Y$ be a finite regular graph of valency, say, $q+1$, $X$ its universal
covering  and
$G$ the automorphism group of $X$.
Then there is a uniform torsion-free lattice $\Ga$ in $G$ such that $Y=\Ga\bs
X$. Recall the zeta function of $Y$,
$$
Z_Y(u)\=\prod_p(1-u^{l(p)})^{-1},
$$
where the product runs over all prime loops in $Y$.
Since $\Ga$ is the fundamental group of $Y$ the set of prime loops is in
bijection to the set of prime conjugacy classes $[\ga_0]$ in $\Ga$.
Hence,
$$
Z_Y(u)\= Z_\Ga(u)\= \prod_{[\ga_0]}(1-u^{l(\ga_0)})^{-1}.
$$
We compute the logarithmic derivative of $Z_Y$,
\begin{eqnarray*}
\frac{Z'_Y}{Z_Y}(u) &=&
\sum_{[\ga_0]}\frac{l(\ga_0)u^{l(\ga_0)-1}}{1-u^{l(\ga_0)}}\\
&=& \sum_{[\ga_0]}l(\ga_0)\sum_{n=1}^\infty \frac{u^{l(\ga_0)n}}{u}\\
&=& \sum_{[\ga]}\frac{u^{l(\ga)}}u l(\ga_0).
\end{eqnarray*}
Here the last sum runs over all nontrivial conjugacy classes in $\Ga$ and
$\ga_0$ is the underlying prime to $\ga$.
The last expression equals the geometric side of the Lefschetz formula for
the test function
$$
\ph(a)\df \frac{u^{l(a)}}u,
$$
which for $|u|<\frac 1q$ satisfies the condition of the Lefschetz formula.
So in that range the function $Z'_Y/Z_Y$ equals
$$
\sum_{a\in A^-}\frac {u^{l(a)}}u q^{l(a)}\,\tr(a|R_S).
$$
For every group homomorphism $\psi$ from $A$ to $\C^\times$ there is a unique number
$\la\in\C^\times$ such that 
$$
\psi(a)\= \la^{v(a)}.
$$
For $\la\in\C^\times$ let $R_S(\la)$ be the \emph{generalized $\la$-eigenspace}, ie,
$R_S(\la)$ is the largest $A$-submodule on which $(a-\la^{v(a)})$ acts
nilpotently for every $a\in A$. Then
$$
R_S\=\bigoplus_{\la}R_S(\la).
$$
Let $m_\la$ denote the dimension of $R_S(\la)$. Then
\begin{eqnarray*}
\frac{Z'_Y}{Z_Y}(u) &=& \sum_{\la} m_\la \sum_{a\in
A^-} \frac{u^{l(a)}}u q^{l(a)}\la^{-l(a)}\\
&=& \sum_{\la} m_\la \sum_{n=1}^\infty
\frac{(uq/\la)^{n}}u\\
&=& \sum_{\la} m_\la \,
\frac{q}{\la-uq}\\
&=& -\sum_{\la} m_\la \,
\frac{1}{u-\la/q}
\end{eqnarray*}
We have proved the following theorem.

\begin{theorem}
The function $Z_Y(u)$ extends to a rational function on $\C$ with $Z_Y(u)$
being a polynomial.
The pole-order of $Z_Y$ at $u=\la/q$, $\la\in\C^\times$, equals
$$
m_\la\=\dim R_S(\la).
$$
\end{theorem}

\section{The Patterson conjecture}
For $\la\in\C^\times$ let $\delta_\la$ be the character of $P=AS$ given by
$$
\delta_\la(as)\df \la^{v(a)}.
$$
Let $I_\la$ denote the representation of $G$ induced from $\delta_\la$. Then $I_\la$ is the
representation on the space of locally constant functions  on $G$ that satisfy
$$
f(px)\= \delta_\la(p)f(x)
$$
for $p\in P$ and $x\in G$. The representation is given by
$$
I_\la(y)f(x)\= f(xy).
$$
The representation $I_\la$ is  admissible. 

\begin{theorem}\label{7.1}
The pole-order of $Z_Y(u)$ at $u=\la$ is
$$
\dim H^0(\Ga,I_{1/\la q}^{-\infty})
$$
if $\la_0\ne \pm \frac 1{\sqrt q}$ and
$$
\dim H^0(\Ga,\hat I_{\pm 1/\sqrt q}^{-\infty})
$$
if $\la=\pm\frac 1{\sqrt q}$, where $\hat I_{\pm 1/\sqrt q}^{-\infty}$ is a certain
self-extension of
$I_{\pm 1/\sqrt q}^{-\infty}$.
\end{theorem}

\begin{corollary}
If $\la\notin \{ \pm \frac 1{\sqrt q}, \pm \frac 1q,\pm 1\}$, then the pole
order of
$Z_Y(u)$ at $u=\la$ equals
$$
\dim H^0(\Ga, I_\la^{-\infty}).
$$
\end{corollary}

\prf
So far we know that the pole-order of $Z_Y(\la)$ equals $\dim
R_S(\la q)=\dim C^\infty(\Ga\bs G)_S(\la q)$.

\begin{lemma}\label{7.2}
Let $\pi$ be an irreducible admissible representation of $G$. If $\la\ne \pm\sqrt q$, then
$A$ acts semisimply on $\pi_S(\la)$. In any case the length of the $A$-module $\pi_S(\la)$
is at most $2$.
\end{lemma}

\prf
Let $\C_\la$ denote the one dimensional $P$-module $\C$ on which $P$ acts by $\delta_\la$.
For $\al\in\Hom(\pi,I_\la)$ let $\al(1)\colon \pi\ra\C$ be given by
$$
\al(1)(v)\= \al(v)(1).
$$
For $s\in S$ we get $\al(1)\pi(s)v=\al(1)v$ and so $\al(1)$ factors over $\pi_S$. In this
way we get an isomorphism
$$
\Hom_G(\pi,I_\la)\ \cong\ \Hom_A(\pi_S,\C_\la).
$$
So, if $\pi_S$ is not zero, then $\pi$, being irreducible, will inject into one $I_\la$
and so $\pi_S$ will inject into $I_{\la,S}$. So it suffices to show the assertion for
$\pi=I_\la$.

The map $\psi\colon \pi=I_\la\ra\C_\la$ given by $f\mapsto f(1)$ is a $P$-homomorphism.
The group $A$ acts on the image by $\delta_\la$ and this defines an irreducible quotient
of $\pi_S$.
We will show that the kernel of $\psi$ gives an irreducible subrepresentation of $\pi_S$
on which $A$ acts by $\delta_{q/\la}$.

Recall the Bruhat decomposition $G=P\cup PwN$.
Let $f$ be in the kernel of $\psi$. Then $f$ vanishes in a neighbourhood of $1$ and thus
the function $s\mapsto f(ws)$ is compactly supported on $S$, so the integral
$\int_Sf(ws)\, ds$ converges.
This integral defines an isomorphism $\ker(\psi)_S\ra\C$.
The group $A$ acts as
\begin{eqnarray*}
\int_Sf(wsa)\, ds &=& q^{v(a)}\int_S f(was)\, ds\\
&=& q^{v(a)}\int_S f(a^{-1}ws)\, ds\\
&=& (q/\la)^{v(a)}\int_S f(ws)\, ds.
\end{eqnarray*}
The lemma follows.
\qed

To be able to apply Lemma \ref{7.2} we have to show that $C^\infty(\Ga\bs G)$ is the
direct sum
$$
C^\infty(\Ga\bs G)\=\bigoplus_{\pi\in\hat G_{adm}} N_\Ga(\pi)\,\pi
$$
of irreducible representations with finite multiplicities.
This follows readily from the corresponding assertion for $L^2(\Ga\bs G)$.
So suppose $\la\ne\pm\frac 1{\sqrt q}$. Then by semisimplicity,
$$
C^\infty(\Ga\bs G)_S(\la q)\= H^0\left(A,C^\infty(\Ga\bs
G)_S\otimes\C_{1/\la q}\right).
$$
Since $A\cong \Z$ we get
\begin{eqnarray*}
\dim H^0\left(A,C^\infty(\Ga\bs G)_S\otimes\C_{1/\la q}\right)
&=& \dim \Hom_A\left(C^\infty(\Ga\bs G)_S,\C_{1/\la q}\right)\\
&=&
\dim \Hom_{G}\left(C^\infty(\Ga\bs G),I_{1/\la q}\right)\\
&=& \dim H^0\left(\Ga,
I_{1/\la q}^{-\infty}\right)
\end{eqnarray*}
In the last equation we have used the duality theorem.
This proves the first assertion of the theorem. For the second replace $A$ with
$A_2=v^{-1}(2\Z)$ and $P$ by the finite index subgroup $A_2S$ and follow the proof above.

For the Corollary use the functional equation
$$
Z_Y\left( \frac 1{qu}\right) \= \left(\frac{1-u^2}{q^2u^2-1}\right)^{r_1-r_0}
q^{2r_1-r_0}\,u^{2r_1}\,Z_Y(u),
$$
which is proved in (\cite{Bass} Cor. 3.10).
Here $r_0$ and $r_1$ are the number of vertices and edges respectively.
\qed

Finally, we give the sheaf-theoretic version of Theorem \ref{7.1}.
For $\la\in\C^\times$, $\la\ne \pm 1/\sqrt q$, let
$$
\CF_\la\df \Ga\bs (X\times I_{1/\la q}^{-\infty}),
$$
where $\Ga$ acts diagonally. If $\la=\pm 1/\sqrt q$ then replace $I$ by $\hat I$
accordingly.
The projection onto the first factor $\CF_\la\ra \Ga\bs X=Y$
makes $\CF_\la$ a locally constant sheaf over $Y$.
Since for  a locally constant sheaf the sheaf-cohomology coincides with the group
cohomology of the fundamental group, we get the following corollary.

\begin{corollary}
The pole order of $Z_Y(u)$ at $u=\la$ equals the dimension of the space of global
sections of the sheaf $\CF_\la$.
\end{corollary}

\newpage
\appendix

\section{Zeta functions for constructible c-sheaves}
In this appendix we suggest a generalization of the Ihara/Bass zeta function. 

A well know generalization is the twist with a finite dimensional representation of the
fundamental group $\Ga$. The latter amounts to the same as the choice of a locally
constant sheaf of vector spaces over the graph. Here we show how the theory can be
extended to a certain class of constructible sheaves. 

A constructible sheaf of a stratified space $X=X_n\supset X_{n-1}\supset\cdots\supset X_0$
is a sheaf which is locally constant on each stratum $X_j\setminus X_{j-1}$.
A graph $X$, viewed as a one dimensional CW-complex, has a natural two-step stratification
$X=X_1\supset X_0$, where $X_0$ is the set of vertices, or the $0$-skeleton. 
Let $\CF$ be a constructible sheaf of abelian groups on $X$.

For each edge $s$ let $x_e$ be its midpoint and let $\CF_e$ denote the stalk at $x_e$.
Let $v$ be an endpoint of $e$. Let $\CF_v$ be the stalk at $v$.
Recall that 
$$
\CF_v=\lim_{\stackrel{\ra}{U\ni v}}\CF(U),
$$
where the limit runs over all open neighbourhoods of $v$. 
Let $s\in \CF_v$. Then there is an open neighbourhood $U$ of $v$ and a section
$s_U\in\CF(U)$ such that $s$ is the class of $s_U$ in $\CF_v$.
Since $U$ is open there is a point $y\in U\cap e$, where $e$ here denotes the open edge.
So $s_U$ also defines a point in the stalk $\CF_y$ over $y$. Since $\CF$ is constant on
the contractible space $e$ there is a canonical isomorphism $\CF_y\cong\CF_e$, so $s_U$
defines an element of $\CF_e$.
Since $\CF$ is constant on $e$ this element does not depend on the choices of $y$ or
$s_U$, but only on $s$. Thus we get a map
$$
\ph_v^e\colon \CF_v\ra\CF_e.
$$
It is easy to see that the sheaf $\CF$ can be recovered up to isomorphism from the stalks
$\CF_v$, $\CF_e$ and the maps $\ph_v^e$.
So a constructible sheaf on a graph is given by the following data
\begin{itemize}
\item
an abelian group $\CF_v$ for every vertex $v$,
\item
an abelian group $\CF_e$ for every edge $e$, and
\item
a group homomorphism $\ph_v^e\colon \CF_v\ra\CF_e$ whenever $v$ is an endpoint of $e$.
\end{itemize}

A \emph{cosheaf} is the dual construction to a sheaf, so a constructible cosheaf is given
abelian groups $\CF_v$, $\CF_e$ as above and group homomorphisms
$$
\psi_e^v\colon \CF_e\ra\CF_v
$$
whenever $v$ is an endpoint of $e$.
Instead of sheaves abelian groups one can also consider sheaves of vector spaces, general
groups, etc.

An example for a cosheaf of groups on a graph $X$ is given as follows.
Let $G$ be the automorphism group of $X$. For each vertex $v$ let $G_v$ be its stabilizer
in $G$. For each edge $e$ let $G_e$ be its pointwise stabilizer in $G$. If $v$ is an
endpoint of $e$ then $G_e$ is a subgroup of $G_v$, so the injection $\psi_e^v\colon G_e\ra
G_v$ defines a constructible cosheaf of groups on $X$.

Next, let $(\rho, V)$ be a representation of $G$. For each vertex $v$ let $\CV_v$ be the
space of fixed points $\CV_v=V^{G_v}$. Likewise for edges.
If $v$ is an endpoint of $e$ then $\CV_v$ is a subspace of $\CV_e$, so the inclusion
$\ph_v^e\colon \CV_v\hookrightarrow \CV_e$ defines a sheaf $\CV$ of vector spaces over $X$.

Assume the representation $\rho$ to be smooth. Then we also have a cosheaf structure
on $\CV$ given by
$$
\begin{array}{cccc}\displaystyle
\psi_e^v\colon & \CV_e & \ra & \CV_v\\
&&&\\
& v&\mapsto & \displaystyle\frac 1{\vol(G_v)} \int_{G_v} \rho(x)v\, dx.
\end{array}
$$
In this case we have $\psi_e^v\ph_v^e = \Id_{\CV_v}$.
This motivates the following definition.

By a \emph{c-sheaf} we mean a constructible sheaf $\CV$ of finite dimensional complex
vector spaces that also carries the structure of a constructible cosheaf such that
$\psi_e^v\ph_v^e = \Id_{\CV_v}$ for every edge $e$ and every endpoint $v$ of $e$.

Let $\CV$ be a c-sheaf on $X$. Let $u,v$ be two adjacent vertices and $e$ their connecting
edge. define $T_u^v=\psi_e^v\ph_u^e$, then $T_u^v$ maps $\CV_u$ linearly to $\CV_v$.
For a path $p=(v_0,\dots,v_n)$ let $T_p\colon\CV_{v_0}\ra\CV_{v_n}$ be given by
$$
T_p\= T_{v_{n-1}}^{v_n}\cdots T_{v_1}^{v_2} T_{v_0}^{v_1}.
$$
Let $c$ le a loop in $X$ given by the path $(v_0,\dots,v_n)$.
Consider the operator
$$
T_c\df T_{v_{n-1}}^{v_n}\cdots T_{v_0}^{v_1}
$$
on $\CV_0$. Note that this operator depends on the choice of the underlying path. However,
for every $j\in\N$ the expression $\tr\left(T_c^j\right)$ does not depend on the choice of
a path.
Therefore also
$$
\det(1-u^{l(c)} T_c)\=\exp\left(-\sum_{j=1}^\infty \frac{u^{jl(c)}}j \,\tr T_c^j\right)
$$
does only depend on the loop $c$.

\begin{theorem}
Let $X$ be a finite graph and $\CV$ a c-sheaf on $X$. Then the infinite product
$$
Z_\CV(u)\df \prod_{p}\det\left(1-u^{l(p)}T_p\right)^{-1},
$$
which runs over all prime loops in $X$, converges for $u$ small and extends to a
rational function in $u$. The inverse $Z_\CV(u)^{-1}$ is a polynomial.
\end{theorem}

\prf
First note that
\begin{eqnarray*}
\tr\left( T_{v_{j+n}}^{v_{j+n+1}}\cdots
T_{v_j}^{v_{j+1}}\right) &=& \tr\left(
\psi_{e_{j+1}}^{v_{j+n+1}}\ph_{v_{j+n}}^{e_{j+1}}
\cdots\psi_{e_j}^{v_{j+1}}\ph_{v_j}^{e_j}\right)\\
&=& \tr \left( \ph_{v_j}^{e_j}
\psi_{e_{j+1}}^{v_{j+n+1}}\ph_{v_{j+n}}^{e_{j+1}}
\cdots\psi_{e_j}^{v_{j+1}}\right)
\end{eqnarray*}
Next consider the set $OE=OE(X)$ of oriented edges. For each oriented edge $e$ let
$e^{-1}$ denote the same edge with reversed orientation and let
$\bar e$ be the underlying simple edge. Define $\CV_e$ to be a copy of $\CV_{\bar e}$ and
let
$$
T\colon \bigoplus_{e\in OE} \CV_e\ra \bigoplus_{e\in OE} \CV_e
$$
be defined by
$$
T(s)\=\sum_{e'} \ph_v^{e'}\psi_e^v(s),
$$
where $s\in\CV_e$ and the sum runs over all oriented edges $e'\ne e^{-1}$ such that the
target point of $e$ equals the start point of $e'$.
It emerges that
$$
\tr T^n\=\sum_{c: l(c)=n} l(c_0)\,\tr T_c.
$$
Therefore,
\begin{eqnarray*}
\det(1-uT) &=&\exp\left( -\sum_{n=1}^\infty\frac{u^n}n\tr T^n\right)\\
&=&\exp\left( -\sum_{n=1}^\infty\frac{u^n}n\sum_{l(c)=n}l(c_0)\,\tr T_c\right)\\
&=&\exp\left( -\sum_{c_0}\sum_{m=1}^\infty \frac{u^{ml(c_0)}}m \,\tr T_{c_0}^m\right)\\
&=& \prod_{c_0}\left( 1-u^{l(c_0)}T_{c_0}\right).
\end{eqnarray*} 
This finishes the proof.
\qed

\section{Graphs with infinite ends}
Linear algebraic groups over local fields of positive characteristic may contain
arithmetic subgroups $\Ga$ which are of finite covolume but not cocompact. The quotient
$\Ga\bs X$ of the Bruhat-Tits building $X$ then is infinite. The finite covolume condition
amounts to the following. For each vertex $v$ of $X$ let $\Ga_v$ be its stabilizer in
$\Ga$. Then $\Ga_v$ is a finite group. It follows that
$$
\sum_{v} \frac 1{|\Ga_v|}\ <\ \infty.
$$
If the building is a tree then the graph $\Ga\bs X$ has the structure described below
\cite{Serre}.

Let $Y$ be a connected graph. Since it does not affect the zeta function we may assume that
$Y$ has no finite ends, ie, no vertices of valency one.
Next we assume the $Y$ decomposes as
$$
Y\= Y_c \cup S_1\cup \dots\cup S_N,
$$
where
$Y_c$ is finite (compact core) and each $S_j$ is a \emph{cusp sector}, ie, the vertices in
$S_j$ are $v_0,v_1,v_2.\dots$, each $v_j$ for $j\ge 1$ has valency $2$, being adjacent
exactly to $v_{j-1}$ and $v_{j+1}$. The intersection of $Y_c$ and $S_j$ is $\{ v_0\}$.

This decomposition is not uniquely determined and the zeta function given below will
depend on the choice.
One can, however, make the choice canonical by insisting, for example, that each
end-vertex $v_0$ of a cusp sector has valency greater that $2$. 
We will now define the zeta function for such a graph in analogy with the Selberg zeta
function for noncompact arithmetic quotients of the upper half plane.

For each vertex $v$ let $d(v)$ denote the distance of $v$ to $Y_c$. 
We say that a path $p=(v_0,\dots,v_n)$ is \emph{essentially backtracking-free} if,
whenever $v_{j-1}=v_{j+1}$, then $d(v_{j-1})<d(v_j)$.
This means that the only backtracking happens in the cusp sectors where one walks in and
turns back once to walk out again.
Next a closed path is \emph{essentially tail-less} if it either is indeed tail-less or if
$d(v_0)>d(v_1)$, ie, the path starts in a cusp sector.
An \emph{S-loop} is an equivalence class of essentially backtracking-free,
essentially tail-less, closed paths under index shift as before.
An S-loop is called \emph{prime} if it is not a power of a shorter one. 

As before let $OE$ denote the set of oriented edges.
We choose a weight $w \colon OE\ra (0,\infty)$ such that
$$
\sum_{e\in OE} w(e)\ <\ \infty.
$$
This is the finite volume condition.
For each S-loop $c$ given by the oriented edges $(e_1,\dots,e_n)$ let
$$
w(c)\df w(e_1)\cdots w(e_n).
$$
Then $w(c^j)=w(c)^j$ for every $j\in\N$.

\begin{theorem}
The infinite product over all prime S-loops,
$$
Z_Y(u)\df \prod_{p}\left(1-u^{l(p)}w(p)\right)^{-1}
$$
converges for $u$ small and extends to a meromorphic function on $\C$. The inverse
$Z_Y(u)^{-1}$ is entire.
\end{theorem}

\prf
Let $\l^2(OE,w)$ be the Hilbert space of all functions $f\colon OE\ra\C$ with
$$
\sum_{e\in OE}|f(e)|^2 w(e) \ <\ \infty.
$$
On $\l(OE,w)$ define a linear operator $T$ by
$$
Tf(e)\df \sum_{e'} f(e')\, w(e),
$$
where the sum runs over all oriented edges $e'$ such that the starting edge of $e'$ equals
the target edge of $e$, and $e'\ne e^{-1}$ except in the case when $e$ lies in a cusp
sector and points outward, in which case $e'=e^{-1}$ is allowed.

\begin{lemma}
The operator $T$ is of trace class and for every $j\in\N$ we have
$$
\tr T^j\= \sum_{l(c)=j} l(c_0)\,w(c),
$$
where the sum runs over all S-loops $c$ of length $j$ and $c_0$ is the underlying prime
S-loop to $c$.
\end{lemma}

\prf 
Once we know that $T$ is of trace class the formula for the trace is clear.
To see that $T$ is of trace class let $e$ be  an oriented edge and define
$\ph_e=\sum_{e'}\delta_{e'}$, where $\delta_{e'}$ is the delta function at $e'$ and the sum
runs over all $e'\in OE$ with start($e'$)=target($e$) and $e'\ne e^{-1}$ except if
$e$ is in a cusp sector and pointing outward in which case $e'=e^{-1}$ is allowed. 
Then we have for $f\in \l(OE,w)$,
$$
Tf\=\sum_{e\in OE} \langle f,\ph_e\rangle\, \delta_e.
$$
To show that $T$ is of trace class it therefore suffices to show that
$$
\sum_{e\in OE} \norm{\ph_e}\norm{\delta_e}\ <\ \infty.
$$
Now $\norm{\delta_e}^2=w(e)$ and $\norm{\ph_e}^2=\sum_{e'}w(e')$, from which it emerges
that $\sum_e\norm{\delta_e}^2<\infty$ and $\sum_e\norm{\ph_e}^2<\infty$, which implies the
claim. The lemma follows.
\qed

Let $\la_1,\la_2,\dots$ be the non-zero eigenvalues of $T$, each repeated according to its
multiplicity. Then $\sum_j|\la_j|<\infty$ and so the product
$$
\det(1-uT)\=\prod_j(1-u\la_j)
$$
converges for every $u\in\C$ and defines an entire function with zeros at the numbers
$\la_j^{-1}$, $j\in\N$.
Using the lemma, the same computation as in the last section yields
$$
\det(1-uT)\= Z_Y(u)^{-1}.
$$
\qed

{\small University of Exeter, Mathematics, Exeter EX4
4QE, England\\ a.h.j.deitmar@ex.ac.uk}

\end{document}